# Strong Solutions of the Fuzzy Linear Systems

Şahin Emrah Amrahov [1] and Iman N. Askerzade [1]

**Abstract**   We consider a fuzzy linear system with crisp coefficient matrix and with an arbitrary fuzzy number in parametric form on the right-hand side. It is known that the well-known existence and uniqueness theorem of a strong fuzzy solution is equivalent to the following: The coefficient matrix is the product of a permutation matrix and a diagonal matrix. This means that this theorem can be applicable only for a special form of linear systems, namely, only when the system consists of equations, each of which has exactly one variable. We prove an existence and uniqueness theorem, which can be use on more general systems. The necessary and sufficient conditions of the theorem are dependent on both the coefficient matrix and the right-hand side. This theorem is a generalization of the well-known existence and uniqueness theorem for the strong solution.

## 1   Introduction

Fuzzy linear systems arise in many branches of science and technology such as economics, social sciences, telecommunications, image processing etc. [Hu et al (2000), Zhang et al (2003), Wu et al (2004), Trivedi et al (2005), Zhou et al (2006), Chen et al (2006), Li at al (2007, 2008)].  In this paper, we consider the fuzzy linear systems whose coefficients matrix is crisp and right-hand side column is an arbitrary fuzzy number vector. These systems have studied by many authors such as, Peeva (1992), Friedman et al (1998, 2000), Asady et al (2005), Abbasbandy et al (2007), Ezzati (2008, 2011), Gasilov et al (2009). There are two approaches for solving the system. Most of the authors try to find the solution in the form of vector of fuzzy numbers. But in this approach, there may be examples having solutions which are not vectors of fuzzy numbers, regardless to which method is used known up to this point. The broadest model in this approach is given by Friedman et al (1998). They transformed the system into $2n \times 2n$ crisp system by using embedding method of Cong-Xin and Ming (1991). After applying the method, if the found solution is a vector of fuzzy numbers, then the solution is called strong solution, if it is not, then it is called weak solution. Friedman et al (2000) also studied duality fuzzy linear systems. They proved necessary and sufficient conditions independently from the right-hand side of the system for the existence of a strong fuzzy solution. Recently, Ezzati (2011) proposed a new method for solving fuzzy linear systems whose coefficients matrix is crisp and right-hand side column is an arbitrary fuzzy number vector. He uses embedding method of Cong-Xin and Ming too. He transforms the Fredman et al's (1998) $2n \times 2n$ system into two $n \times n$ systems.

---

[1] Ankara University, Computer Engineering Department

Gasilov et al (2009) suggested a second approach for linear systems, which is a geometric approach based on the properties of linear transformations. Unlike the first approach, in this approach, the solution is considered to be a fuzzy set of vectors. Therefore, the suggested approach can be applied to all linear square systems. Gasilov et al (2009) apply this approach to non-square linear systems and differential equation systems too (2011). Comparison analysis between these two approaches is done by Wierzchon (2010).

In this paper we analyze Friedman et al's existence and uniqueness theorem for strong solution. It is known that these brought up conditions are equivalent to the following: The coefficient matrix can be obtained by multiplying a permutation matrix by a diagonal matrix. This shows that these conditions are applicable only in narrow cases. That is, conditions hold for only when the system consists of equations, each of which has exactly one variable. Geometric proof of the theorem of Friedman et al is handled by Gasilov et al (2009). For obtaining an existence theorem which can be used on more general systems; we searched for conditions dependent on both the coefficient matrix and the right-hand side. The existence and uniqueness theorem of the strong solution proved in this paper is in fact more general and coincides with the theorem of Friedman et al (1998) in special cases.

This paper is organized in 6 sections. In Section 2, we state some basic definitions of fuzzy sets and numbers. In Section 3, we define fuzzy linear systems (FLS) and transform FLS to linear systems of crisp equations as done in Friedman et al (1998, 2000). In this section we also define strong and weak solutions. In Section 4, we prove the necessary and sufficient conditions for existence of unique strong solution of the auxiliary crisp linear system corresponding to FLS. In Section 5, we give analogous conditions for FLS. In Section 6 we conclude the paper with some remarks.

## 2  Preliminaries

As Dubois et al (1978) we define a fuzzy number $\tilde{u}$ in parametric form.

**Definition 1**.  A fuzzy number $\tilde{u}$ in parametric form is a pair $(\underline{u},\overline{u})$ of functions $\underline{u}(r), \overline{u}(r), 0 \leq r \leq 1$, which satisfy the following requirements:

1. $\underline{u}(r)$ is a bounded monotonically increasing left continuous function over [0,1]
2. $\overline{u}(r)$ is a bounded monotonically decreasing left continuous function over [0,1]
3. $\underline{u}(r) \leq \overline{u}(r), 0 \leq r \leq 1$

The set of all these fuzzy numbers is denoted by $E^1$.

A popular type of fuzzy number is triangular numbers $\tilde{u}=(a,c,b)$ with the membership function $\mu(x) = \begin{cases} \dfrac{x-a}{c-a}, a \leq x \leq c \\ \dfrac{x-b}{c-b}, c \leq x \leq b \end{cases}$ where $c \neq a, c \neq b$. For triangular numbers we have $\underline{u}(r) = a + (c-a)r$ and $\overline{u}(r) = b + (c-b)r$.

We can represent a crisp number $a$ by $\overline{u}(r) = \underline{u}(r) = a$, $0 \leq r \leq 1$.

**Definition 2**. For two arbitrary fuzzy numbers $\tilde{u}$ and $\tilde{v}$ the equality $\tilde{u} = \tilde{v}$ means that

$\underline{u}(r) = \underline{v}(r)$

and

$\overline{u}(r) = \overline{v}(r)$ for all $r \in [0,1]$.

**Definition 3**. For two arbitrary fuzzy numbers $\tilde{u}$ and $\tilde{v}$ addition is defined as

$$\tilde{u} + \tilde{v} = ((\underline{u+v})(r), \overline{(u+v)}(r)) = (\underline{u}(r) + \underline{v}(r), \overline{u}(r) + \overline{v}(r))$$

**Definition 4.** For any fuzzy number $\tilde{u}$ and real number $k > 0$ multiplication by positive real number is defined as

$$k\tilde{u} = (k\underline{u}(r), k\overline{u}(r))$$

**Definition 5.** For any fuzzy number $\tilde{u}$ and real number $k < 0$ multiplication by negative real number is defined as

$$k\tilde{u} = (k\overline{u}(r), k\underline{u}(r))$$

### 3 Fuzzy linear systems

**Definition 6.** Let a crisp $n \times n$ coefficient matrix $A = (a_{ij}), 1 \leq i, j \leq n$ and fuzzy numbers $b_i \in E^1, 1 \leq i \leq n$ be given. Then the $n \times n$ system of equations

$$
\begin{aligned}
a_{11}x_1 + a_{12}x_2 + \ldots a_{1n}x_n &= b_1, \\
a_{21}x_1 + a_{22}x_2 + \ldots a_{2n}x_n &= b_2, \\
&\vdots \\
a_{n1}x_1 + a_{n2}x_2 + \ldots a_{nn}x_n &= b_n
\end{aligned}
\quad (1)
$$

is called a fuzzy linear system (FLS).

**Definition 7.** A fuzzy number vector $(x_1, x_2, \ldots x_n)^T$ is called a solution of the fuzzy system if the fuzzy numbers $x_1, x_2, \ldots x_n$ satisfy the system (1) in the sense of definitions 2-5.

Following Friedman et al (1998) we introduce the notations below:

$$x = (\underline{x}_1, \underline{x}_2, \ldots \underline{x}_n, -\overline{x}_1, -\overline{x}_2, \ldots -\overline{x}_n)^T$$

$$b = (\underline{b}_1, \underline{b}_2, \ldots \underline{b}_n, -\overline{b}_1, -\overline{b}_2, \ldots -\overline{b}_n)^T$$

$S = (s_{ij}), 1 \leq i, j \leq 2n$, where $s_{ij}$ are determined as follows:

$$
\begin{aligned}
a_{ij} \geq 0 &\Rightarrow s_{ij} = a_{ij}, \; s_{i+n, j+n} = a_{ij}, \\
a_{ij} < 0 &\Rightarrow s_{i, j+n} = -a_{ij}, \; s_{i+n, j} = -a_{ij}
\end{aligned}
\quad (2)
$$

and any $s_{ij}$ which is not determined by (2) is zero. Using matrix notation we have

$$Sx = b \tag{3}$$

The structure of S implies that $s_{ij} \geq 0$ and that

$$S = \begin{pmatrix} B & C \\ C & B \end{pmatrix} \tag{4}$$

where $B$ contains the positive elements of $A$, $C$ contains the absolute value of the negative elements of $A$ and $A = B - C$. An example in the work of Friedman et al (1998) shows that the matrix $S$ may be singular even if $A$ is nonsingular.

**Theorem1** (Friedman et al (1998)) The matrix S is nonsingular if and only if the matrices $A = B - C$ and $B + C$ are both nonsingular.

The following example shows that the solution of the crisp linear system (3) does not define a fuzzy solution of the system (1) even if S is nonsingular.

**Example**. Consider the following fuzzy system

$$\begin{aligned} x_1 - x_2 &= b_1 \\ x_1 + 2x_2 &= b_2 \end{aligned} \tag{5}$$

where $b_1 = (\underline{b}_1(r), \overline{b}_1(r))$ and $b_2 = (\underline{b}_2(r), \overline{b}_2(r))$.

We will write briefly $b_1 = (\underline{b}_1, \overline{b}_1)$, $b_2 = (\underline{b}_2, \overline{b}_2)$.

The extended matrix is

$$S = \begin{pmatrix} 1 & 0 & 0 & 1 \\ 1 & 2 & 0 & 0 \\ 0 & 1 & 1 & 0 \\ 0 & 0 & 1 & 2 \end{pmatrix} \tag{6}$$

Since $A = \begin{pmatrix} 1 & -1 \\ 1 & 2 \end{pmatrix}$ and $B + C = \begin{pmatrix} 1 & 1 \\ 1 & 2 \end{pmatrix}$ are nonsingular matrices by the Theorem 1, $S$ is nonsingular matrix. Therefore the crisp system (3) has a unique solution. The solution of the system (3) for the matrix (6) is

$$\underline{x}_1 = \underline{b}_1 + \frac{2}{3}(\overline{b}_2 - \overline{b}_1) - \frac{1}{3}(\underline{b}_2 - \underline{b}_1)$$

$$\overline{x}_1 = \overline{b}_1 - \frac{1}{3}(\overline{b}_2 - \overline{b}_1) + \frac{2}{3}(\underline{b}_2 - \underline{b}_1)$$

$$\underline{x}_2 = -\frac{1}{3}(\overline{b}_2 - \overline{b}_1) + \frac{2}{3}(\underline{b}_2 - \underline{b}_1)$$

$$\overline{x}_2 = \frac{2}{3}(\overline{b}_2 - \overline{b}_1) - \frac{1}{3}(\underline{b}_2 - \underline{b}_1)$$

Note that $x_1 = (\underline{x}_1, \overline{x}_1)$ and $x_2 = (\underline{x}_2, \overline{x}_2)$ are not necessarily fuzzy numbers. From the condition $\underline{x}_1 \leq \overline{x}_1$ we have

$$\overline{b}_2 - \underline{b}_2 \leq 2(\overline{b}_1 - \underline{b}_1) \tag{7}$$

From the condition $\underline{x}_2 \leq \overline{x}_2$ we get

$$\overline{b}_1 - \underline{b}_1 \leq \overline{b}_2 - \underline{b}_2 \tag{8}$$

Therefore $x_1 = (\underline{x}_1, \overline{x}_1)$ and $x_2 = (\underline{x}_2, \overline{x}_2)$ are fuzzy numbers if and only if the inequalities (7) and (8) hold. As in Friedman et al (1998) if the solution of the crisp system (3) defines fuzzy solution to the fuzzy system (1) this solution is called strong solution.

**Definition 8**. If $x = (\underline{x}_1, \underline{x}_2, ... \underline{x}_n, -\overline{x}_1, -\overline{x}_2, ... -\overline{x}_n)^T$ is a solution of (3) and for each $1 \leq i \leq n$ the inequalities

$$\underline{x}_i \leq \overline{x}_i$$

hold, then the solution $x = (\underline{x}_1, \underline{x}_2, ...\underline{x}_n, -\overline{x}_1, -\overline{x}_2, ... - \overline{x}_n)^T$ is called a strong solution of the system (3).

**Definition 9**. If $x = (\underline{x}_1, \underline{x}_2, ...\underline{x}_n, -\overline{x}_1, -\overline{x}_2, ... - \overline{x}_n)^T$ is a solution of (3) and for some $i \in [1, n]$ the inequality

$$\underline{x}_i > \overline{x}_i$$

hold, then the solution $x = (\underline{x}_1, \underline{x}_2, ...\underline{x}_n, -\overline{x}_1, -\overline{x}_2, ... - \overline{x}_n)^T$ is called a weak solution of the system (3).

## 4 Necessary and sufficient conditions for the existence of a strong solution

Let us define

$$\underline{b} = (\underline{b}_1, \underline{b}_2, ...\underline{b}_n) \tag{9}$$

and

$$\overline{b} = (\overline{b}_1, \overline{b}_2, ...\overline{b}_n) \tag{10}$$

**Theorem 2**. Let $S = \begin{pmatrix} B & C \\ C & B \end{pmatrix}$ be a nonsingular matrix. The system (3) has a strong solution if and only if

$$(B+C)^{-1}(\underline{b}-\overline{b}) \leq 0 \tag{11}$$

**Proof**. Similar to (9) and (10), we can define:

$$\underline{x} = (\underline{x}_1, \underline{x}_2, ...\underline{x}_n)$$

$$\overline{x} = (\overline{x}_1, \overline{x}_2, ...\overline{x}_n)$$

From the system (3) we obtain:

$$\begin{pmatrix} B & C \\ C & B \end{pmatrix} \begin{pmatrix} \underline{x} \\ -\overline{x} \end{pmatrix} = \begin{pmatrix} \underline{b} \\ -\overline{b} \end{pmatrix}$$

Hence

$$B\underline{x} - C\overline{x} = \underline{b} \tag{12}$$

$$-C\underline{x} + B\overline{x} = \overline{b} \tag{13}$$

From (12) and (13) we have

$$(B+C)\underline{x} - (B+C)\overline{x} = \underline{b} - \overline{b}$$

$$(B+C)(\underline{x} - \overline{x}) = \underline{b} - \overline{b}$$

By the Theorem 1, the matrix $B+C$ is nonsingular. Therefore

$$\underline{x} - \overline{x} = (B+C)^{-1}(\underline{b} - \overline{b}) \tag{14}$$

If the system (3) has a strong solution then, by the Definition 8, we have $\underline{x} - \overline{x} \leq 0$. Hence the inequality (11) holds. Conversely, if the inequality (11) holds, by (14), we have $\underline{x} - \overline{x} \leq 0$

## 5 The necessary and sufficient conditions for the existence of a unique strong solution of FLS

By the Theorems 1 and 2, we have the following result:

**Theorem 3.**

The FLS (1) has a unique strong solution if and only if the following conditions hold:

1) The matrices

$A = B - C$ and $B + C$ are both nonsingular.

2) $(B+C)^{-1}(\underline{b} - \overline{b}) \leq 0$

**Example.** Consider the following fuzzy system again.

$$x_1 - x_2 = b_1$$
$$x_1 + 2x_2 = b_2$$

Here $b_1 = (\underline{b}_1(r), \overline{b}_1(r))$ and $b_2 = (\underline{b}_2(r), \overline{b}_2(r))$.

For this example, we have the matrices

$$A = \begin{pmatrix} 1 & -1 \\ 1 & 2 \end{pmatrix}$$

$$B = \begin{pmatrix} 1 & 0 \\ 1 & 2 \end{pmatrix}$$

$$C = \begin{pmatrix} 0 & 1 \\ 0 & 0 \end{pmatrix}$$

$$B + C = \begin{pmatrix} 1 & 1 \\ 1 & 2 \end{pmatrix}$$

$$(B+C)^{-1} = \begin{pmatrix} 2 & -1 \\ -1 & 1 \end{pmatrix}$$

The matrices $A = \begin{pmatrix} 1 & -1 \\ 1 & 2 \end{pmatrix}$ and $B + C = \begin{pmatrix} 1 & 1 \\ 1 & 2 \end{pmatrix}$ are nonsingular matrices. Applying the Theorem 3, we note that for existence of a unique strong fuzzy solution, the necessary and sufficient condition is

$$\begin{pmatrix} 2 & -1 \\ -1 & 1 \end{pmatrix} \begin{pmatrix} \underline{b_1} - \overline{b_1} \\ \underline{b_2} - \overline{b_2} \end{pmatrix} \leq 0$$

Hence

$$\overline{b_1} - \underline{b_1} \leq \overline{b_2} - \underline{b_2} \leq 2(\overline{b_1} - \underline{b_1}).$$

**Corollary.** The system (1) has a unique fuzzy solution for an arbitrary $b = (b_1, b_2, ... b_n)^T$ if and only if

$(B+C)^{-1}$ is nonnegative, in other words

$(B+C)^{-1}_{ij} \geq 0, 1 \leq i, j \leq n$

This result is equivalent to the Lemma 2 of Friedman et al (2000) which states that the system (1) has a unique fuzzy solution for arbitrary $b = (b_1, b_2, ... b_n)^T$ if and only if $S^{-1}$ is nonnegative. But it turns out that this condition restricts the coefficient matrix $A$ to a very specific case. The entries of the matrices $S$ and $(B+C)$ are nonnegative too. It is a well-known fact from linear algebra [Anton et al (2003)] that if a nonsingular matrix and its inverse are both nonnegative matrices (i.e. matrices with nonnegative entries), then the matrix is a generalized permutation matrix (monomial matrix). Furthermore, a nonsingular matrix is a generalized permutation matrix if and only if it can be written as a product of a nonsingular diagonal matrix and a permutation matrix. Hence we can apply this corollary or Lemma 2 of Friedman et al (2000) only on a small class of FLS.

## 6  Conclusion

Friedman and et al (1998, 2010) proved several necessary and sufficient conditions independent of right-hand side for the existence of a strong solution to the FLS. In this paper, we point out that these conditions are applicable only in certain narrow cases. We have suggested a generalized version of these conditions which are additionally dependent on the right-hand side of the system.